%% file: osid_BMC17082020.tex
\documentclass[11pt]{book}
\usepackage{graphicx}

\catcode`\@=11

 \def\@normalsize{\@setsize\normalsize{13pt}\xipt\@xipt
   \abovedisplayskip 11pt plus3pt minus6pt
   \belowdisplayskip \abovedisplayskip
   \abovedisplayshortskip \z@ plus3pt
   \belowdisplayshortskip 6.6pt plus3.5pt minus3pt}

 \def\small{\@setsize\small{12pt}\xipt\@xipt
   \abovedisplayskip 10pt plus2pt minus5pt
   \belowdisplayskip \abovedisplayskip
   \abovedisplayshortskip \z@ plus3pt
   \belowdisplayshortskip 6pt plus3pt minus3pt
   \def\@listi{\topsep 6pt plus 2pt minus 2pt
     \parsep 3pt plus 2pt minus 1pt
     \itemsep \parsep}}

 \def\footnotesize{\@setsize\footnotesize{10pt}\ixpt\@ixpt
   \abovedisplayskip 8pt plus 2pt minus 4pt
   \belowdisplayskip \abovedisplayskip
   \abovedisplayshortskip \z@ plus 1pt
   \belowdisplayshortskip 4pt plus 2pt minus 2pt
   \def\@listi{\topsep 4pt plus 2pt minus 2pt
      \parsep 2pt plus 1pt minus 1pt
      \itemsep \parsep}}

 \def\scriptsize{\@setsize\scriptsize{9.5pt}\viiipt\@viiipt}
 \def\tiny{\@setsize\tiny{7pt}\vipt\@vipt}
 \def\large{\@setsize\large{14pt}\xiipt\@xiipt}
 \def\Large{\@setsize\Large{18pt}\xivpt\@xivpt}
 \def\LARGE{\@setsize\LARGE{22pt}\xviipt\@xviipt}
 \def\huge{\@setsize\huge{25pt}\xxpt\@xxpt}
 \def\Huge{\@setsize\Huge{30pt}\xxvpt\@xxvpt}

\def\section{\@startsection {section}{1}{\z@}%
{-1.5\baselineskip plus-1pt minus-3pt}{1\baselineskip plus1pt minus2pt}%
{\centering\normalsize\bf}}
\def\subsection{\@startsection{subsection}{2}{\z@}%
{-1\baselineskip plus-1pt minus-2pt}{1\baselineskip plus1pt minus2pt}%
{\normalsize\sc\noindent}}
\def\subsubsection{\@startsection{subsubsection}{3}{\z@}%
{-1\baselineskip plus-1pt minus-2pt}{1sp}{\normalsize\it\noindent}}
\def\paragraph{\@startsection{paragraph}{4}{\z@}%
{1\baselineskip plus1pt minus2pt}{-1em}{\normalsize\it\noindent}}
\let\subparagraph=\paragraph

\def\tableofcontents{\@restonecolfalse\if@twocolumn\@restonecoltrue
\onecolumn\fi\OSIDcont\@starttoc{con}\if@restonecol\twocolumn\fi}

\def\l@section{\@dottedtocline{1}{0em}{.66em}}

\def\thebibliography#1{\section*{{Bibliography}\@mkboth
 {BIBLIOGRAPHY}{BIBLIOGRAPHY}}\footnotesize\rm\list
 {[\arabic{enumi}]}{\settowidth\labelwidth{[#1]}\leftmargin\labelwidth
 \advance\leftmargin\labelsep\usecounter{enumi}}
 \def\newblock{\hskip .11em plus .33em minus -.07em}
 \sloppy\clubpenalty4000\widowpenalty4000
 \sfcode`\.=1000\relax}

\def\ps@myheadings{\let\@mkboth\@gobbletwo
\def\@oddhead{\hfil{\footnotesize\rm\rightmark}\hfil}
\def\@evenhead{\hfil{\footnotesize\rm\leftmark}\hfil}
\def\@oddfoot{\hfil{\footnotesize\sf\artid-\thepage}\hfil}
\def\@evenfoot{\hfil{\footnotesize\sf\artid-\thepage}\hfil}
\def\sectionmark##1{}\def\subsectionmark##1{}}

\def\@copyrighthead{\parbox{127mm}{\footnotesize\rm\ \\[6pt]
Open Systems~\& Information Dynamics\\
Vol.~\Vol, No.~\Number~(\Year)~\artid~(\EndpagE~pages)\\
DOI:\DOInumber\\
\copyright~World Scientific Publishing Company\\
\epsfxsize=4cm
\vskip-\lastskip
\vskip-\baselineskip
\vspace*{-38.5pt}
\noindent\hfill\epsfbox{wlogo.eps}}}

\def\artid{0000001}
\def\Year{2008}        %
\def\Vol{15}           
\newcounter{paPer}     %
\def\EndpagE{\expandafter\pageref{\the\value{paPer}OpSy}}

\def\ps@osiD{\let\@mkboth\@gobbletwo
\def\@oddhead{\@copyrighthead}
  \def\@oddfoot{\hfil{\footnotesize\sf\artid-\thepage}\hfil}
  \def\@evenhead{}\let\@evenfoot\@oddfoot}

\def\cite{\@ifnextchar [{\@tempswatrue\@Rcitex}{\@tempswafalse\@Rcitex[]}}

\def\@Rcitex[#1]#2{\if@filesw\immediate\write\@auxout{\string\citation{#2}}\fi
  \def\@citea{}\@cite{\@for\@citeb:=#2\do
    {\@citea\def\@citea{,\penalty\@m\,}\@ifundefined
       {b@\the\value{paPer}R\@citeb}{{\bf ?}\@warning
       {Citation `\@citeb' on page \thepage \space undefined}}%
\hbox{\csname b@\the\value{paPer}R\@citeb\endcsname}}}{#1}}

\long\def\@caption#1[#2]#3{\par\addcontentsline{\csname
  ext@#1\endcsname}{#1}{\protect\numberline{\csname
  the#1\endcsname}{\ignorespaces #2}}\begingroup
    \@parboxrestore
    \small                                        
    \@makecaption{\csname fnum@#1\endcsname}{\ignorespaces #3}\par
  \endgroup}

\newtoks\@stequation

\def\subequations{\refstepcounter{equation}%
\edef\@savedequation{\the\c@equation}%
\@stequation=\expandafter{\theequation}
\edef\@savedtheequation{\the\@stequation}
\edef\oldtheequation{\theequation}%
\setcounter{equation}{0}%
\def\theequation{\oldtheequation\alph{equation}}}%

\def\endsubequations{%
\setcounter{equation}{\@savedequation}%
\@stequation=\expandafter{\@savedtheequation}%
\edef\theequation{\the\@stequation}\global\@ignoretrue}

\catcode`\@=12

\let\Rlabel=\label
\let\Rbibitem=\bibitem
\let\Rref=\ref
\let\Rpageref=\pageref
\def\label#1{\expandafter\Rlabel{\the\value{paPer}R#1}}
\def\bibitem#1{\expandafter\Rbibitem{\the\value{paPer}R#1}}
\def\ref#1{\expandafter\Rref{\the\value{paPer}R#1}}
\def\pageref#1{\expandafter\Rpageref{\the\value{paPer}R#1}}

\def\thesection{\arabic{section}.}

\def\YYMm{\rule{0ex}{4em}}
\newtoks\TITsi
\newtoks\TITsii

\def\title#1{\def\TITs{\LARGE{\raggedright\noindent\YYMm #1%
\vskip8pt\par}}}

\def\author#1{\autMM{#1}\def\LHD{#1}}
\def\and{{\rm\lowercase{and}}}

\def\autMM#1{\TITsii={\vskip10pt\par\normalsize\rm\noindent #1\par}%
\TITsi=\expandafter{\TITs}\edef\TITs{\the\TITsi\the\TITsii}}

\def\address#1{\TITsii={\vskip6pt\par\footnotesize\sl\noindent #1\par}%
\TITsi=\expandafter{\TITs}%
\edef\TITs{\the\TITsi\the\TITsii}}

\def\received#1{\TITsii={\vskip10pt\par\small\rm\noindent(Received: #1)\par}%
\TITsi=\expandafter{\TITs}\edef\TITs{\the\TITsi\the\TITsii}}

\def\headtitle#1{\def\RHD{#1}}

\def\abst{{\bf Abstract.}}
\def\abstract#1{\TITs
       \vskip15pt\par\noindent
       {\footnotesize{\abst~} #1\vskip3pt\par}
       \markright{\RHD}
       \markboth{\LHD}{\RHD}}

\def\OSIDcont{\cleardoublepage\thispagestyle{empty}
       \markright{}\markboth{}{}
       \normalsize\rm
       \hspace*{\fill}{\large\rm
         Contents of the Volume \Volume, Number \Number}\hspace*{\fill}
       \par\vspace{1.5em}
       \par\noindent}

\def\endpaper{\expandafter\label{\the\value{paPer}OpSy}}

\def\emptyset{\mathchoice{\mbox{\normalsize\rm\O}}{\mbox{\normalsize\rm\O}}
{\mbox{\scriptsize\rm\O}}{\mbox{\tiny\rm\O}}}

\def\1{{\mathchoice{\rm 1\mskip-4mu l}{\rm 1\mskip-4mu l}%
{\rm 1\mskip-4.5mu l}{\rm 1\mskip-5mu l}}}

\def\varkappa{\mbox{\bBB\char 123}}

\def\longhookrightarrow{\lhook\joinrel\relbar\joinrel\rightarrow}

\def\longhookUp{\lower6pt\hbox{\rotatebox{90}{$\longhookrightarrow$}}}

\newtheorem{thm}{\rm THEOREM}
\newtheorem{cor}{\rm COROLLARY}
\newtheorem{lem}{\rm LEMMA}

\newtheorem{prop}{\rm PROPOSITION}

\newtheorem{defn}{\rm DEFINITION}
\newtheorem{rem}{\it Remark}

\def\LD{\Lambda}
\def\theequation{\thesection\arabic{equation}}

\def\Myskip{\setlength{\baselineskip}{13pt}}

\def\text#1{\quad\mbox{\rm  #1 }\quad}

\def\DOInumber{}

\input xy
\xyoption{all}

\InputIfFileExists{psfig.sty}{\typeout{^^Jpsfig.sty inputed...ok}}{\typeout{^^JWarning: psfig.sty could not be found.^^J}}
\def\mathfrak{\mathbf}
\begin{document}

\input{00Ch_rot.tex}

\end{document}

%% file: 00Ch_rot.tex
\newcommand{\Mn}{M_n(\mathbb{C})}
\newcommand{\Mk}{M_k(\mathbb{C})}
\newcommand{\id}{\mbox{id}}
\newcommand{\ot}{{\,\otimes\,}}
\newcommand{{\Cd}}{{\mathbb{C}^d}}
\newcommand{\sbsigma}{{\mbox{\scriptsize \boldmath $\sigma$}}}
\newcommand{\sbalpha}{{\mbox{\scriptsize \boldmath $\alpha$}}}
\newcommand{\sbbeta}{{\mbox{\scriptsize \boldmath $\beta$}}}
\newcommand{\bsigma}{{\mbox{\boldmath $\sigma$}}}
\newcommand{\balpha}{{\mbox{\boldmath $\alpha$}}}
\newcommand{\bbeta}{{\mbox{\boldmath $\beta$}}}
\newcommand{\bmu}{{\mbox{\boldmath $\mu$}}}
\newcommand{\bnu}{{\mbox{\boldmath $\nu$}}}
\newcommand{\ba}{{\mbox{\boldmath $a$}}}
\newcommand{\bb}{{\mbox{\boldmath $b$}}}
\newcommand{\sba}{{\mbox{\scriptsize \boldmath $a$}}}
\newcommand{\MD}{\mathfrak{D}}
\newcommand{\sbb}{{\mbox{\scriptsize \boldmath $b$}}}
\newcommand{\sbmu}{{\mbox{\scriptsize \boldmath $\mu$}}}
\newcommand{\sbnu}{{\mbox{\scriptsize \boldmath $\nu$}}}
\def\oper{{\mathchoice{\rm 1\mskip-4mu l}{\rm 1\mskip-4mu l}%
{\rm 1\mskip-4.5mu l}{\rm 1\mskip-5mu l}}}
\def\<{\langle}
\def\>{\rangle}
\def\theequation{\thesection\arabic{equation}}
\def\mathbb{\mathbf}

\title{Block Markov Chains on Trees }
\author{Abdessatar Souissi}
\address{$^1$ Department of Accounting, College of Business Management\\
Qassim University, Ar Rass, Saudi Arabia \\
$^2$ Preparatory Institute for Scientific and Technical Studies La Marsa,\\
 Carthage University, Tunisia}

\headtitle{Block Markov Chains on Trees}

\abstract{We introduce block Markov chains (BMCs) indexed by an infinite rooted tree.
It turns out that BMCs define a new class of tree-indexed Markovian  processes. We clarify the structure of BMCs in connection with Markov chains (MCs) and Markov random fields (MRFs). Mainly,  show that probability measures which  are  BMCs for every root are indeed
  Markov chains (MCs) and yet they form a strict subclass of Markov random fields (MRFs) on the considered tree. Conversely, a class of MCs which are BMCs is characterized. Furthermore, we establish that in the one-dimensional case the class of  BMCs coincides
with MCs. However, a slight perturbation of the one-dimensional lattice leads to us to an example of BMCs which are not MCs appear.}

\Myskip



\section{Introduction}\label{sect_intro}
 Markov random fields (MRFs) on lattice  have become standard tools in several branches  of science and technology including computer science, machine learning,  graphical models, statistical physics. Namely, MRFs are known to provide pertinent models for interacting particles systems in statistical mechanics. \\
We notice that MRFs were introduced by Dobrushin  in \cite{D68} for the multi-dimensional integer lattice, and developed then on trees   \cite{Spi75}, \cite{Spa75}, \cite{Za83}. QMFs consist   multi-dimensional extensions of Markov chains  \cite{Norris}  but with a deeper Markovian structure.  In, fact even in the one dimensional case MRFs were shown to be distinct from MCs \cite{ChGMMP2013}.

 MRFs play a crucial  role in many areas such as  computer science, image recognition, graphical models, psychology and in an increasing number of biological and neurological models. The reader is referred to \cite{Lau98}, \cite{Coa}, \cite{Gly} and the references cited therein for further applications.\\
In the present paper we introduce the notion of block Markov chains  indexed by the vertex-set of a  rooted tree $T= (V, E)$. The definition of this notion is quite natural. Since in  the one-dimensional case $V= \mathbb N_0$ with distinguished vertex (root) $"o=0"$,   a  Markov chain
 $(Z_n)_{n\in \mathbb N}$ with (finite) state space $\Xi$ is  defined through the well known Markov property
$$
\mathbb P[Z_{n+1}= \xi_{n+1} \, \mid  \, Z_{n}= \xi_{n}, \cdots, Z_{0}= \xi_{0} ]  =  \mathbb P[Z_{n+1}= \xi_{n+1} \, \mid  \, Z_{n}= \xi_{n}].
$$
The above property can be reformulated by means of the joined probability measure $\mu$  on  $\Xi^V$ of the process $(Z_u)_{u\in V}$ as follows
\begin{equation}\label{1st_eq_BMC}
\mu[ \xi(.) \, \,\hbox{on} \, \,   S(x) \, \mid  \,   \xi(.) \, \,  \hbox{on} \, \,  V\setminus T^{'}(x)) ]
 = \mu[\xi(.) \, \,\hbox{on} \, \, S_(x) \, \mid  \, \, \xi(x)]
\end{equation}
where $S(x) =\{x+1\}$ is the set of successors of the site $x\in V$ and $T'(x) = \{x+1, x+2, \cdots\}$ it the set of successive descendants
  of the vertex  $x$ w.r.t. the considered root $o$. We emphasize a suitable natural generalization of the sets $S(x)$ and  $T(x)$
 for general trees. Roughly speaking, a BMC is a probability  measure on $\Omega:= \Xi^V$ satisfying the Markov property (\ref{1st_eq_BMC}) for a fixed root.\\
 The main purpose of this paper is to clarify the structure of BMCs in connection with MCs and MRFs. Mainly, we show that a probability measure which is BMC for every root $o\in V$ is a MC in the sense of \cite{Za83}. The correlation functions of BMCs are different from those of MCs and MRFs. Consequently, their Markov structure are also different.  Namely, it turns out that some conditional independence conditions are necessary on a MC on the considered tree to be BMC.\\
On the other hand we show that in the one-dimensional case, the   notions of MCs and QMCs coincide. This coincidence Makes BMCs strictly a sub-class of MRFs in the one dimensional case. However, we emphasize that a slight modification of the one-dimensional lattice leads to a counter-examples that confirms  the  huge difference between  MCs and BMCs over multi-dimensional trees.\\
We notice that the natural hierarchical structure of rooted trees, due to the absence of loops,  plays  a crucial role in the mere definition of BMCs. Therefore, the results are no longer available on general graphs.  We forecast that BMCs will  play a crucial role in connection with  Gibbs measures on trees and their associated phenomena of phase transitions (see  \cite{Rozikov}, \cite{Minlos2014} and \cite{Geor}). Namely, phenomena of phase transitions were associated with  interesting p-adic models such as the Potts model and the Ising–Vannimenus model \cite{MuAk13}, \cite{MuSaKh15}. 
In fact, a work under preparation is dedicated to the clarification of a bridge between BMCs and some p-adic models.\\
In \cite{MS19}, \cite{MS20} we clarified the structure of quantum Markov states on  a quasi-local algebra $\mathcal{A}$ trees in terms of classical Markovian measure and  Gibbs measures on the spectrum of a maximal abelian subalgebra. We stress that this classical Markovian measure is indeed a BMC. This will makes a new bridge between classical and quantum Markov fields. \\
Let us mention the outlines of the paper.  Section \ref{sect_tree} is devoted to some notions and notions on rooted trees.
In section \ref{sect_mf}, we recall the basic definition of  MC and MRF on graphs.
Section \ref{sect_BMC} is devoted to definition  of BMCs as far as its correlation functions. Section \ref{Sect_main} is dedicated to results related to the connection of BMCs with MCs and MRFs on trees.
In section \ref{sect_1D} we deal with the one-dimensional case for which the vertex set is the classical 1D integer
lattice $\mathbb Z$ occupied with its natural tree structure.   In section \ref{sect_Cexple} we develop a counter-example for a BMC which is not a MC.

 \section{Rooted trees}\label{sect_tree}
 Recall that \cite{Ru_Ly16} a tree is a connected graph with no cycles, i.e. a connected graph which becomes disconnected when each one of its edges is removed.\\
Let be given  an infinite tree  $T = (V,E)$. First, we fix  any vertex $o=x_0\in V$ as a "root".  Recall that two vertices $x$ and $y$ are said to be  {\it nearest neighbors}and we denote $x\sim y$ if they are joined through an edge (i.e. $<x,y>\in E$). A list of the vertices $ x\sim x_1\sim \dots \sim x_{d-1}\sim y$ is called a {\it
path} from the site $x$ to the site $y$. The distance $d(x,y)$ on the  tree  is the length of the shortest path from $x$ to $y$. \\
For  $x\in V$, its \textit{direct successors} (children)  is defined  by
\begin{equation}\label{S(x)def}
S^o(x) :  = \left\{y\in V \, \,  : \, \,  x\sim y \, \, \hbox{and} \, \, d(y,o) > d(x,o) \right\}
\end{equation}
and its $k^{th}$ successors w.r.t. the root $o$  is defined by induction as follows
$$
S_1^o(x) := S^o(x);
$$
$$
  S^o_{k+1}(x) = S^o(S^o_k(x)),\, \,   \forall k\ge 1.
$$
The "future" w.r.t. the vertex $x$ is defined by:
\begin{equation}\label{S_1n_S_infty}
  S^o_{[m,n]}(x)= \bigcup_{k=m}^{n}S^o_{k}(x); \quad    T_o(x) = \bigcup_{k\ge1}S^o_{k}(x); \quad T_o^{'}(x) = T(x)\setminus \{x\}.
\end{equation}
 Note that in the homogeneous case, for which $|S_o(x)| = k$ is constant, the graph $T$ is  the semi-infinite Cayley tree $\Gamma^k_+$ of order $k$.
Namely, for $k=1$, the graph is reduced to the one-dimensional integer lattice $\mathbb Z$.\\
Consider the map $r$ from $V$ into itself characterized by
$$
r(o) =o,
$$
$$
  r(y) =x \quad     \hbox{if}\quad  y\in S^o(x)
$$
Let $x\in V$. If $n= d(x, o)$ then

\begin{equation}\label{roots_x}
o= r^{n}(x) = x_0 \sim r^{n-1}(x) \sim \cdots \sim r(x) \sim r^0(x) = x
\end{equation}
is the minimal edge-path joining  the root $o$ to the vertex $x$, where $r^k =\underbrace{ r\circ \cdots\circ r}_{k \, \,  times}$. \\
The set
\begin{equation}\label{R(x)}
R^o(x):= \{r(x), r^2(x), \cdots, r^{n}(x)=o \}
\end{equation}
 represents the "past"  of the vertex $x$ for the root $o$\\
 The set of {\it nearest-neighbors vertices} of $x$ is given as follows:
\begin{equation}\label{N_x}
  N_x = \{y \in V \, \,  : \, \,  x\sim y \}
\end{equation}
It is clear that $N_x = \{r(x)\}\cup S^o(x)$.\

 In the sequel, the tree $T$ is assumed to be locally finite, i.e. $|N_x|<\infty$ for each $x\in V$, in this case the integer $d_x:= |N_x|$ is called \textit{degree} of $x$.\\
The tree can be regarded as growing (upward) away from its fixed root $o$. Each vertex  $x\in V$  then
has branches leading to its "children", which are represented here by $Sô(x)$ and $T^{'}_o(x)$. With the possibility of leaves, that is, vertices $x$ without children i.e.  $S(x)= \emptyset$.

\section{Some Reminders on Markov fields}\label{sect_mf}
Let $\Xi = \{1,\cdots, q \}$. By a stochastic process we mean a family of random variables  $(Z_{u})_{u\in V}$  defined on a probability space
$(\Omega, \mathcal F, \mathbb P)$ and valued in a finite set $\Xi:= \{1,2,\cdots, q\}$. The process $(Z_{u})_{u\in V}$ is  defined through its joined probability measure $\mu$ on the Borel  space $( \Xi^V, \mathcal B)$ where  $\mathcal B_V$ is the cylindrical $\sigma$-algebra, which is generated by the cylinder sets of the following form
 \begin{equation}\label{cylinder}
 C( a_x, \; x\in \Lambda ) = \left\{ \xi\in \Xi^V\, \,  : \, \,  \xi(x) =a_x, \, \,  \forall x\in \Lambda\right\}
 \end{equation}
  where $\Lambda\subset V$ finite and $(a_x)_{x\in \Lambda}\in \Xi^{|\Lambda|}$. For the sake of shortness we denote $\Omega$ instead of $\Omega_V$  and  $\mathcal B$ instead of $\mathcal{B}_V$. For $\Lambda\subset V$, we denote $\Omega_\Lambda = \Xi^\Lambda$.
Recall that
\begin{equation}\label{Z_u_mu_coressp}
  \mu\left[\xi(.) \,\, \hbox{on} \, \,  \Lambda \right] = \mathbb P\left[ \bigcap_{u\in \Lambda}( Z_u = \xi(u ))\right]
\end{equation}
where $\xi\in\Omega_\Lambda$.\\
The conditional probability is defined as follows
\begin{equation}\label{mu_Z mid_coressp}
   \mu\left[\xi(.) \,\, \hbox{on} \, \,  \Lambda \, \, \mid \, \, \xi(.) \,\, \hbox{on} \, \,  \Lambda^{'} \right] =
 \frac{ \mu\left[\xi(.)  \,\, \hbox{on} \, \, \Lambda\cup \Lambda^{'} \right]}{\mu\left[\xi(.)  \,\, \hbox{on} \, \,  \Lambda^{'} \right]}
\end{equation}
where $\Lambda, \Lambda^{'}\subseteq V$ and $\xi \in \Xi^{V}$ such that
$$ \mu\left[\xi(.)  \,\, \hbox{on} \, \,  \Lambda^{'} \right]>0.$$
Denoting
\begin{equation}\label{F(u)_F(lb)}
\mathcal F_u := \sigma(Z_u) \, \,   ; \, \,  \mathcal F_\Lambda =\sigma \left( Z_u\, \, ; \, \, u\in \Lambda\right)
\end{equation}
the $\sigma$-algebra generated by $Z_u$ and $(Z_v, \, v\in\Lambda),$ respectively.

\begin{defn}\cite{D68}\label{MF_def}
A probability measure $\mu$ on $(\Omega, \mathcal{B})$ is said to be   Markov  random  field (MRF)  if it takes strictly positive values on finite cylinder sets of the form (\ref{cylinder}) and such that for every $\xi \in \Omega$
\begin{equation}\label{l_MF_eq}
  \mu\left[\xi(u) \, \mid \xi(.) \;  \hbox{on}   \; V\setminus\{u\}  \right] = \mu\left[ \xi(u) \,  \mid  \,  \xi(. ) \,  \hbox{on} \,  N_u \right].
\end{equation}
The set of Markov random fields over $T$ will be denoted by $\mathcal{MF}(T)$.
\end{defn}
The conditional probabilities  (\ref{l_MF_eq}) are assumed to be invariant under graph isomorphism.

\begin{defn}\label{def_Mc}\cite{Spa75}A probability measure $\mu$ on $(\Omega, \mathcal{B})$  is said to  be  Markov chain (MC) over the tree $T=(V,E)$ if for each subtree $T'= (V', E')$   the restriction of $\mu$ on  the measurable space $(\Omega_{V'}, \mathcal B_{V^{'}})$ defines a Markov random field. i.e.
\begin{equation}\label{MC_eq}
  \mu[\xi(x) \, \, \mid \, \, \xi(.) \, \, \hbox{on} \, \, V'\setminus\{x\} \, ] = \mu[\xi(x) \, \, \mid \, \xi(.) \, \, \hbox{on} \, N_x \cap V^{'} \, ]
\end{equation}
for all $x\in V'$ and all $\xi \in \Omega_{V^{'}}$.
The set of Markov chains over $T$ will be denoted by $\mathcal{MC}(T)$.
\end{defn}
\begin{rem}
The class $\mathcal{MC}(T)$ is clearly  included in $\mathcal{MF}(T)$. Conversely, in \cite{Za83} it was proven that if the tail $\sigma$-field is trivial
then the considered Markov field is indeed a MC.
\end{rem}

\section{Structure of Block Markov chains on trees}\label{sect_BMC}

In what follows, a root $o$ for the tree $T= (V,E)$ is fixed. For each $n\in\mathbb{N}$, we denote  $\LD_n := S_{n}^o(o)$   the set of vertices whose
 distance to the root $o$ equals $n$.
Let  $\Lambda_{n]} = S_{[0,n]}^o(o) = \bigcup_{k=0}^{n}\LD_k$. For the sake of shortness, when  confusion seems  impossible we will
  use the notations  $ S(x),\,  T(x), T^{'}$  and $r (x) $  instead of $S_o(x), \, T_o(x), T_o^{'}(x)$ and $r_o(x)$, respectively.

Let us set a random enumeration for elements of $\LD_n$ as follows
\begin{eqnarray*}
{\LD_n}:=\left(x^{(1)}_{\LD_n},x^{(2)}_{\LD_n},\cdots,x^{(|\LD_n|)}_{\LD_n}\right)
\end{eqnarray*}
where $|\LD_n|$ denotes the cardinality of $\LD_n$.

\begin{defn}\label{defbMC}
A measure $\mu$ on $(\Omega, \mathcal{B})$ is called  o-block Markov chain (o-BMC) if  it satisfies
\begin{equation}\label{Mp_x_S_o(x)}
\mu\left[\xi(.)\, \,  \hbox{on} \; S(x)\, \,  \big| \xi(.)  \, \,  \hbox{on}\; V\setminus T^{'}(x)  \right] = \mu\left[\xi(.)\, \,  \hbox{on}\, \,  S(x)\; \big| \xi(x)  \right]
\end{equation}
 for all $x\in V$ and $\xi\in \Omega$. The equation (\ref{Mp_x_S_o(x)}) will be referred as \textit{block Markov property}.
 The set of $o$-block Markov chains over the tree $T$ will be denoted $o-\mathcal{BMC}(T)$.
\end{defn}

In \cite{Spa75}  a triplet of $\sigma$-algebras $(\mathcal F_1, \, \mathcal F_2, \, \mathcal F_3)$ such that
\begin{equation}\label{Mt}
\mathbb P(A \, \mid \, \mathcal F_1 \vee \mathcal F_2) = \mathbb P(A \, \mid \, \mathcal F_2), \quad \forall A\in \mathcal{F}_3
\end{equation}
 was referred as \textit{Markov triple}. In these notations  (\ref{Mp_x_S_o(x)}) means that $( \mathcal F_{S(x)},
\, \mathcal F_{V\setminus T(x)},\,  \mathcal{F}_x )$ is  a Markov triple.\\
\begin{rem}
 The word "block" in Definition. \ref{defbMC} comes from the conditioning w.r.t. the   $\sigma$-algebra $\mathcal F_{V\setminus T(x)}$ rather then the
   $\sigma$-algebra $\mathcal F_{R(x)}$, while this latter  represents the past of the vertex $x$ w.r.t. the root $o$.
   \end{rem}
 The following  elementary formula for conditional probabilities will be used frequently in the sequel.
\begin{equation}\label{condi_prob_formula}
\mathbb P(A\cap B \mid C ) =  \mathbb P(A\mid B\cap C) \mathbb P(B\mid C).
\end{equation}
Let  $\mu$ is an $o$-BMC. According to (\ref{condi_prob_formula}), we have
\begin{eqnarray*}
\mu[\xi(.) \, \hbox{on}\, \Lambda_{n]} ]& = & \mu[\xi(.) \, \hbox{on}\, \LD_n \, \mid \,\xi(.) \, \hbox{on}\,  \Lambda_{n-1]}]\times\mu[\xi(.) \, \hbox{on}\, \Lambda_{n-1]} ]\\
&=& \mu[\xi(.) \, \hbox{on}\, \Lambda_{0} ] \prod_{k=0}^{n-1} \mu[\xi(.) \, \hbox{on}\, \LD_{k+1} \, \mid \,\xi(.) \, \hbox{on}\,  \Lambda_{k]}].
\end{eqnarray*}
For $k=1,\cdots, n-1$, the same reason as above implies that
$$
\mu[\xi(.) \, \hbox{on}\, \LD_{k+1} \, \mid \,\xi(.) \, \hbox{on}\,  \Lambda_{k]}] =
\prod_{j=1}^{|\LD_k|}\mu\bigl[\xi(.) \, \hbox{on} \, S(x_{\LD_k}^{(j)}) \, \mid \, \xi(.) \, \hbox{on} \, \Lambda_{n-1]}\cup\bigcup_{i=j+1}^{|\LD_k|}S(x_{\LD_k}^{(i)})\bigr].
$$
Since $x_{\LD_k}^{(i)}\in \Lambda_{n-1]}\cup\bigcup_{i=j+1}^{|\LD_k|}S(x_{\LD_k}^{(i)})\subset V \setminus T^{'}(x_{\LD_k}^{(i)})$ then the block Markov property
 (\ref{Mp_x_S_o(x)}) leads to
$$
\mu\bigl[\xi(.) \, \hbox{on} \, S(x_{\LD_k}^{(j)}) \, \mid \, \xi(.) \, \hbox{on} \, \Lambda_{n-1]}\cup\bigcup_{i=j+1}^{|\LD_k|}S(x_{\LD_k}^{(i)})\bigr]
 = \mu\bigl[\xi(.) \, \hbox{on} \, S(x_{\LD_k}^{(j)}) \, \mid \, \xi(x_{\LD_k}^{(j)}) \bigr].
$$
Therefore
\begin{equation}\label{bMC_onLb_n}
  \mu[ \xi(.) \, \hbox{on} \, \Lambda_{n]}] = \mu[\xi(o)]\prod_{k=0}^{n-1}\prod_{x\in \LD_k}\mu\bigl[\xi(.) \, \hbox{on} \, S(x) \, \mid \, \xi(x) \bigr].
\end{equation}
\begin{rem}
The BMC $\mu$  is  characterized by the initial distribution $\mu_o$ on $\Omega_{\{o\}}$ together with the family of transition probabilities $\mu\bigl[\xi(.) \, \hbox{on} \, S(x) \, \mid \, \xi(x) \bigr]$. The $ d\times (d^{|S(x)|})$ "stochastic" matrices
$\Pi_{x, S(x)} = \left(\mu[\xi^{'}(.) \, \hbox{on} \, S(x) \, \mid \, \xi(x)]\right)_{\xi^{'}\in \Xi^{S(x)}, \xi\in \Xi^{\{x\}}}$
 are clearly inhomogeneous. This lets the measure $\mu$ a multi-dimensional markovian process which is inhomogeneous  both in space and time.
\end{rem}

The following theorem extends the local Markov property (\ref{Mp_x_S_o(x)}) to a global one, which concerns the conditional independence
of the $\sigma$-algebras  $\mathcal F_{ T(x) }$ and $\mathcal F_{ V\setminus T^{'}(x)} $ given $\mathcal F_{x} $.\\
\begin{thm}\label{Sinfty_midxtheorem} Let  $\mu$ be a block Markov chain on $(\Omega, \mathcal{B})$. Then
\begin{equation}\label{S_inft(x)_mid_x}
  \mu\left[\xi(.) \, \hbox{on} \, \, T^{'}(x)\, \big|\,  \xi(.) \, \hbox{on} \, \,    V\setminus T^{'}(x) \right] =  \mu\left[\xi(.) \, \hbox{on} \, \,  T^{'}(x)\, \big|\, \xi(x)\right]
\end{equation}
For all $\xi\in \Omega$ and all $x\in V$.
\end{thm}

\textit{Proof}.   If $T^{'}(x) = \emptyset$ then (\ref{S_inft(x)_mid_x}) holds true.\\
  We will proceed by induction on $S_{\left[1,n\right]}(x): = \bigcup_{k=1}^{n}S_{k}(x)$. One has

   $$\mu\left[\xi(.) \, \,  \hbox{on} \, \,  S_{\left[1,n+1\right]}(x) \, \,  \big| \, \,  \xi(.) \, \,  \hbox{on} \, \,    V\setminus T^{'}(x)\right] $$
 $$  =   \mu\left[\xi(.) \, \,  \hbox{on} \, \,  S_{n+1}(x)\, \,   \big| \, \,  \xi(.) \, \,  \hbox{on} \, \,  S_{\left[1,n\right]}(x)\cup  V\setminus T^{'}(x) \right]$$
 $$\times \mu\left[\xi(.) \, \,  \hbox{on} \, \,  S_{\left[1,n\right]}(x)\, \,  \big| \, \,   \xi(.) \, \,  \hbox{on} \, \,    V\setminus T^{'}(x)\right].$$
Denoting $S_{n}(x) =\{x_1^{(n)}, \cdots, x^{(n)}_{|S_{n}(x)|}\}$, one has
$$
 \mu\left[\xi(.) \, \,  \hbox{on} \, \,  S_{n+1}(x)\, \,   \big| \, \,  \xi(.) \, \,  \hbox{on} \, \,  S_{\left[1,n\right]}(x)\cup   V\setminus T^{'}(x) \right]
 $$
 $$
 = \prod_{i=1}^{|S_{n}(x) |} \mu\left[\xi(.) \, \,   \hbox{on} \, \,  S(x_i^{(n)}) \, \,   \big| \, \,  \xi(.) \, \,  \hbox{on} \, \,  ( \bigcup_{k=i+1}^{n} S(x_k^{(n)}))\cup S_{\left[1,n\right]}(x)\cup  V\setminus T^{'}(x) \right].
$$
 From (\ref{Mp_x_S_o(x)}), one gets
$$
\mu\left[\xi(.) \, \,   \hbox{on} \, \,  S(x_i^{(n)}) \, \,   \big| \, \,  \xi(.) \, \,  \hbox{on} \, \,  ( \bigcup_{k=i+1}^{n} S(x_k^{(n)}))\cup S_{\left[1,n\right]}(x)\cup  V\setminus T^{'}(x) \right]
$$
$$
= \mu\left[\xi(.)   \,  \, \hbox{on} \, \,  S(x_i^{(n)}) \, \,   \big| \, \,  \xi(x_i^{(n)}) \right].
$$
Thus
$$
 \mu\left[\xi(.) \, \,  \hbox{on} \, \,  S_{n+1}(x)\, \,   \big| \, \,  \xi(.) \, \,  \hbox{on} \, \,  S_{\left[1,n\right]}(x)\cup   V\setminus T^{'}(x) \right]
 $$
 $$
 = \prod_{i=1}^{|S^{(n)}(x) |} \mu\left[\xi(.) \, \,   \hbox{on} \, \,  S(x_i^{(n)}) \, \,   \big| \, \,  \xi(x_i^{(n)}) \right].
$$
Using the same argument as above , we obtain
\begin{equation}\label{S_n+1mid _S_n}
  \mu\left[\xi(.) \, \,  \hbox{on} \, \,  S_{n+1}(x)\, \,   \big| \, \,  \xi(.) \, \,  \hbox{on} \, \,  S_{\left[1,n\right]}(x) \right]
 \end{equation}
 $$
 =\prod_{i=1}^{|S_{n}(x) |} \mu\left[\xi(.) \, \,   \hbox{on} \, \,  S(x_i^{(n)}) \, \,   \big| \, \,  \xi(x_i^{(n)}) \right].
 $$
 On the other hand, the induction's hypthesis leads to
 $$
 \mu\left[\xi(.) \, \,  \hbox{on} \, \,  S_{[1,n]}(x)\, \,  \big| \, \,   \xi(.) \, \,  \hbox{on} \, \,    V\setminus T^{'}(x)\right] =
  \mu\left[\xi(.) \, \,  \hbox{on} \, \,  S_{[1,n]}(x)\, \,  \big| \, \,   \xi(x) \right].
 $$
 Therefore
 $$
  \mu\left[\xi(.) \, \,  \hbox{on} \, \,  S_{[1,n+1]}(x)\, \,   \big| \, \,  \xi(.) \, \,  \hbox{on} \, \,    V\setminus T^{'}(x) \right]
  $$
  $$
  = \mu\left[\xi(.) \, \,  \hbox{on} \, \,  S_{n+1}(x)\, \,   \big| \, \,  \xi(.) \, \,  \hbox{on} \, \,  S_{\left[1,n\right]}(x) \right]\times
  \mu\left[\xi(.) \, \,  \hbox{on} \, \,  S_{[1,n]}(x)\, \,  \big| \, \,   \xi(x) \right]
 $$
 $$
 =  \mu\left[\xi(.) \, \,  \hbox{on} \, \,  S_{[1, n+1]}(x)\, \,    \big| \, \,   \xi(x) \right].
 $$
Finally, one finds
  \begin{eqnarray*}
 && \mu\left[\xi(.) \, \,  \hbox{on} \, \,  T^{'}(x)\, \,   \big| \, \,  \xi(.) \, \,  \hbox{on} \, \,    V\setminus T^{'}(x)
 \right]\\& =&
 \lim_{n\to[1,\infty)}
  \mu\left[\xi(.) \, \,  \hbox{on} \, \,  S_{[1,n+1]}(x)\, \,   \big| \, \,  \xi(.) \, \,  \hbox{on} \, \,    V\setminus T^{'}(x) \right];\\
   & =& \lim_{n\to[1,\infty)}
  \mu\left[\xi(.) \, \,  \hbox{on} \, \,  S_{[1,n+1]}(x)\, \,   \big| \, \,  \xi(x) \right];\\
   &  = &\mu\left[\xi(.) \, \,  \hbox{on} \, \,  T^{'}(x)\, \,   \big| \, \,  \xi(x) \right].
 \end{eqnarray*}

 \begin{cor}\label{cor_Ld-mid-x}
In the  notations of Theorem \ref{Sinfty_midxtheorem},  if $\Lambda\subseteq T^{'}(x)$ then
 \begin{equation}\label{Lbd8x}
  \mu\left[\xi(.) \, \hbox{on} \, \,   \Lambda \, \,   \big|\,  \xi(.) \, \hbox{on} \, \,    V\setminus T^{'}(x) \right] =  \mu\left[\xi(.) \, \hbox{on} \, \,   \Lambda \, \big| \,  \xi(x)\right]
\end{equation}
for all $\xi\in \Omega$.
 \end{cor}
 \textit{Proof.}
From  Theorem \ref{Sinfty_midxtheorem}, for each $\xi'\in \Omega_{T^{'}(x)\setminus\Lambda}$
 $$
  \mu\left[\xi(.) \, \,   \hbox{on} \, \,   \Lambda ,\, \,  \xi^{'}(.)\, \,   \hbox{on} \, \,   T^{'}(x)\setminus\Lambda \, \,  \big| \, \,   \xi(.) \, \hbox{on} \, \,    V\setminus T^{'}(x)\right]
   $$
   $$
   =  \mu\left[\xi(.) \, \,   \hbox{on} \, \,   \Lambda, \, \,  \xi^{'}(.)\, \,   \hbox{on} \, \,   T^{'}(x)\setminus \Lambda\, \,  \big| \, \,   \xi(x)\right].
 $$
 Summing up on $ {\xi'\in \Omega_{T^{'}(x)\setminus\Lambda}}$, one finds  (\ref{Lbd8x}).  \\
 The following result proposes a multi-dimensional  analogue of the Chapmann-Kolmogorov equation.
\begin{thm}\label{theorem_chap_Kol}
Let $\mu$ be a BMC on $(\Omega, \mathcal{B})$. Then for $x\in V$ and  $m,n \in \mathbb N$ one has
\begin{equation}\label{Cap_kol}
\mu\left[ \xi(.) \, \, \hbox{on} \, \, S_{n+m}(x) \, \, \big| \, \, \xi(x) \right]
\end{equation}
$$
= \sum_{\xi^{'}\in \Xi^{S_{n}(x)}}\mu\left[ \xi(.) \, \, \hbox{on} \, \, S_{n+m}(x) \, \, \big|\xi^{'}(.) \, \, \hbox{on} \, \, S_{n}(x) \, \, \right]\times \mu\left[ \xi^{'}(.) \, \, \hbox{on} \, \, S_{n}(x) \, \, \big|\xi(x) \, \, \right]
$$
for all $\xi\in \Omega$.
\end{thm}
\textit{Proof.}
For each $\xi^{'}\in \Omega_{S_n (x)}$, using the same reason as in (\ref{S_n+1mid _S_n}), we get
$$
\mu\left[ \xi(.) \, \, \hbox{on} \, \, S_{n+m}(x) \, \, \big|\xi^{'}(.) \, \, \hbox{on} \, \, S_n(x) \, \, \right]
$$
$$
=\mu\left[ \xi(.) \, \, \hbox{on} \, \, S_{n+m}(x) \, \, \big|\xi^{'}(.) \, \, \hbox{on} \, \, S_{n}(x) \, , \,  \xi(x) \,  \right].
$$
Then
$$
\mu\left[ \xi(.) \, \, \hbox{on} \, \, S_{n+m}(x) \, \, \big|\xi^{'}(.) \, \, \hbox{on} \, \, S_n(x) \, \, \right]\times \mu\left[ \xi^{'}(.) \, \, \hbox{on} \, \, S_n(x) \, \, \big|\, \, \xi(x) \, \, \right]
$$
$$
= \left[ \xi(.) \, \, \hbox{on} \, \, S_{n+m}(x) \, , \, \xi^{'}(.) \, \, \hbox{on} \, \, S_n(x) \, \,  \big|\, \, \xi(x) \, \, \right].
$$
Summing up, one gets (\ref{Cap_kol}).
\section{Connection with MCs and MRFs}\label{Sect_main}
 \begin{lem}\label{lem_subtree}
Let $x\in V$. If  $\Lambda$ is a  subset of $S(x)$ then the subgraph of the tree $T= (V,E)$, whose set of vertices is $\Lambda\cup (V\setminus T^{'}(x))$ is itself a tree.
 \end{lem}
 \textit{Proof.}
 First, we see that  if $y\in  T^{'}(x)$  then   $T^{'}(y)\subseteq  T^{'}(x)$. This implies that for
 each $y\in V\setminus T^{'}(x) $ is connected, the set of its roots $\{r^k(y), \, k =0, \cdots\}$ ( defined in  (\ref{roots_x}))
  is disjoint of the set  $T^{'}(x)$. Then the path $y\sim r(y)\sim\cdots \sim o$ is in  $V\setminus T^{'}(x) $.
Therefore, the  subgraph whose vertex set $V\setminus T^{'}(x)$ is connected. Since every element of $S(x)$ is joined to $x$, we conclude that the subgraph $(\Lambda\cup (V\setminus T^{'}(x)), \sim)$ is connected. Taking into account that the  fact that every connected subgraph of a tree is a subtree, the proof is complete.\\
\begin{thm}\label{theoremMCbMC} Let $\mu$ be a Markov chain on $\Omega$. Then for each $x\in V$ the following property holds true.
 \begin{equation}\label{MC_bMC_eq}
\mu\left[\xi(.)\, \,  \hbox{on}\; S(x)\, \,  \big| \, \, \xi(.)  \, \,  \hbox{on}\; V\setminus T^{'}(x)  \right] = \prod_{y\in S(x)}\mu\left[\xi(y) \, \, \mid \, \, \xi(x)\right].
\end{equation}
 If in addition, the $\sigma$-algebras  $(\mathcal F_y)_{ y\in S(x)}$ are conditionally independent given $\mathcal F_x$ then $\mu$ is  an o-BMC.
\end{thm}
\textit{Proof.}
 First let us write $ S(x):= \{y_1, y_2, \cdots, y_{|S(x)|} \}$.
 According to (\ref{condi_prob_formula}), we have
 $$
 \mu\left[\xi(.)\, \,  \hbox{on}\; S(x)\, \,  \big| \, \, \xi(.)  \, \,  \hbox{on}\; V\setminus T^{'}(x)  \right]
  $$
  $$
  = \prod_{k=1}^{|S(x)|}\mu\left[\xi(y_k) \, \, \mid \, \, \xi(.) \, \,  \hbox{on}\; (V\setminus T^{'}(x))\cup \{y_{k+1}, \cdots, y_{n}\}\right].
 $$
By Lemma \ref{lem_subtree} the subgraph of $T$ whose set of vertices is $V^{'}:= (V\setminus T^{'}(x))\cup\{y_{k+1}, \cdots, y_{n}\}$ is a tree. Then
 \begin{eqnarray*}
 && \prod_{k=1}^{|S(x)|}\mu\left[\xi(y_k) \, \, \mid \, \, \xi(.) \, \,  \hbox{on}\; (V\setminus T^{'}(x))\cup \{y_{k+1}, \cdots, y_{n}\}\right]\\
& =& \prod_{k=1}^{|S(x)|}\mu\left[\xi(y_k) \, \, \mid \, \, \xi(.) \, \,  \hbox{on}\; (V\setminus\{y_k\})\cap V^{'}\right]\\
& =& \prod_{k=1}^{|S(x)|}\mu\left[\xi(y_k) \, \, \mid \, \, \xi(.) \, \,  \hbox{on}\; N_{y_k}\cap V^{'}\right].
 \end{eqnarray*}
 where  the last equality derives from the fact that $\mu$ is a Markov chain in the sense of Definition. \ref{def_Mc}.
 Since $ N_{y_k}\cap V^{'}= \{x\}$, we get  (\ref{MC_bMC_eq}).
 For the second part of the proof, the conditional independence of $\mathcal F_y:= \sigma(Z_y), y\in S(x)$ leads to
 $$
  \prod_{k=1}^{|S(x)|}\mu\left[\xi(y) \, \, \mid \, \, \xi(x) \right] =  \mu\left[\xi(.) \, \,  \hbox{on}\; S(x) \, \mid \, \, \xi(x) \right].
 $$
 Hence, (\ref{MC_bMC_eq}) leads to (\ref{S_inft(x)_mid_x}). Therefore $\mu$ is a o-block Markov chain, for any root $o\in V$. This achieves the proof.\\
\begin{lem}\label{lem_Bmc_mf}
If $\mu$ is an $o$-BMC on $ (\Omega, \mathcal{B})$  and $x\in V$ then
\begin{equation}\label{x_mid xC}
\mu[\xi(x)\, \,  \mid \, \,  \xi \, \, \hbox{on} \, \Lambda ] = \mu[\xi(x)\, \,  \mid \, \xi(.) \, on  \, \{r(x)\}\cup ( T^{'}(x)\cap \Lambda)]
\end{equation}
for all $\Lambda \subseteq V\setminus\{x\}$ containing $r(x)$.
\end{lem}
\textit{Proof.} Since $x\notin \Lambda$, then  according to (\ref{S_inft(x)_mid_x})  one gets
\begin{eqnarray*}
\mu\left[ \xi(x) \, \,  \big| \, \,  \xi(.) \, \,  \hbox{on} \, \, \Lambda\right]
&=& \frac{\mu\left[  \xi(x)\, ; \, \,\xi(.) \, \,  \hbox{on} \, \,    \Lambda\cap T^{'}(x)\, ; \; \xi(.) \, \,  \hbox{on} \, \,   \Lambda \setminus T^{'}(x) \right] }
{\mu\left[ \xi(.) \, \,  \hbox{on} \, \,    \Lambda\cap T^{'}(x)\, ; \, \,  \xi(.) \, \,  \hbox{on} \, \,    \Lambda \setminus  T^{'}(x)\right]}, \\
&=& \frac{\mu\left[\xi(.) \, \,  \hbox{on} \, \,    \Lambda\cap T^{'}(x) \, \,  \big| \, \,     \xi(x)\right] \mu\left[ \xi(x) \, \,  \big| \, \,   \xi(.) \, \,  \Lambda \setminus T^{'}(x)  \right] }
{\mu \left[ \xi(.) \, \,  \hbox{on} \, \,   \Lambda\cap T^{'}(x) \, \,  \big| \, \,   \xi(.) \, \,  \hbox{on} \, \,   \Lambda \setminus T^{'}(x)  \right]}\\
&=& \frac{\mu\left[\xi(.) \, \,  \hbox{on} \, \,    \Lambda\cap T^{'}(x) \, \,  \big| \, \,     \xi(x)\right] \mu\left[ \xi(x) \, \,  \big| \, \,   \xi(r(x))\right] }
{\mu \left[ \xi(.) \, \,  \hbox{on} \, \,    \Lambda\cap T^{'}(x) \, \,  \big| \, \,   \xi(r(x)) \right]}.
\end{eqnarray*}
Again from  (\ref{S_inft(x)_mid_x}), we have
$$
\mu \left[ \xi(.) \, \,  \hbox{on} \, \,    \Lambda\cap T^{'}(x) \, \,  \big| \, \,   \xi(x) \right] = \mu \left[ \xi(.) \, \,  \hbox{on} \, \,    \Lambda \cap T^{'}(x) \, \,  \big| \, \,   \xi(x) \, \,  ; \, \,  \xi(r(x)) \right].
$$
This leads to
\begin{eqnarray*}
\mu\left[ \xi(x) \, \,  \big| \, \,  \xi(.) \, \,  \hbox{on} \, \,  \Lambda\right]
&=&\frac{\mu\left[\xi(x) \, \,  ; \, \,   \xi(.) \, \,  \hbox{on} \, \,    \Lambda\cap T^{'}(x) \, \,  , \, \,   \xi(r(x))\right] }
{\mu \left[ \xi(.) \, \,  \hbox{on} \, \,   \Lambda\cap T^{'}(x) \, \, ; \, \,   \xi(r(x)) \right]}\\
&=&\mu\left[ \xi(x) \, \,  \big| \, \,  \xi(.) \, \,  \hbox{on}  \, \,  \{r(x)\} \cup \Lambda\cap T^{'}(x)\right].
\end{eqnarray*}
This completes the proof.\\
\begin{rem}
Notice that Definition.\ref{def_Mc} extends the notion of Markov chain introduced in \cite{Spi75} and \cite{Spa75} into inhomogeneous trees and for inhomogeneous transition probabilities. It was shown \cite{Spi75} that the class of homogenous  Markov chain is strictly included in the class of Markov random fields. In the inhomogeneous we have the following
\end{rem}
\begin{thm}\label{theorem_bMC_MC}
Let $\mu$ be a probability measure on $(\Omega, \mathcal{B})$. If $\mu$ is an $o$-BMC for each $o\in V$ then it is a MC.
\end{thm}
\textit{Proof.}  Consider a subtree $T^{'} = (V^{'}, E^{'})$ of $T$. Let $x\in V^{'}$. If $V^{'}\cap N_x = \emptyset$ then $V^{'} = \{x\}$
and (\ref{MC_eq}) is trivial. Otherwise,  let us denote $N_x\cap V' = \{y_1, y_2, \cdots, y_d\}$ with $d= |N_x\cap V^{'}|$.
Remark that if $o=y\in N_x\cap V^{'}$  then $r_o(x) = y$
 and $N_x\cap V^{'} \setminus\{y\} \subseteq S(x)\cap V^{'} \subseteq T^{'}_o(x)\cap V^{'}$. As $\mu$ is an $y_{1}$-BMC,
 by Lemma \ref{lem_Bmc_mf} we have
$$
\mu[\xi(x) \, \, \mid  \, \, \xi(.) \, \hbox{on} \, \,V^{'}\setminus  \{x\} ] =
 \mu[ \xi(x) \, \, \mid \, \, \xi(.) \, \hbox{on} \, \, \{y_1\}\cup ( T_{y_1}^{'}(x)\cap V^{'}) ].
$$
Since $\mu$ is an $y_2$-BMC   by Lemma \ref{lem_Bmc_mf}  , we have
\begin{eqnarray*}
&&\mu[ \xi(x) \, \, \mid \, \, \xi(.) \, \hbox{on} \, \, \{y_1\}\cup (T_{y_1}^{'}(x)\cap V^{'}) ]\\
&=& \mu[ \xi(x) \, \, \mid \, \, \xi(.) \, \hbox{on} \, \, \{y_2\}\cup \bigl(\{y_1\}\cup T_{x_1}^{'}(x) \cap T_{y_2}^{'}(x)\cap V^{'}\bigr) ]\\
&=& \mu[ \xi(x) \, \, \mid \, \, \xi(.) \, \hbox{on} \, \, \{y_1, y_2\}\cup \bigl( T_{y_1}^{'}(x) \cap T_{y_2}^{'}(x)\cap V^{'}\bigr) ].
\end{eqnarray*}
Iterating this procedure, we get
\begin{eqnarray*}
\mu[\xi(x) \, \, \mid  \, \, \xi(.) \, \hbox{on} \, \,V^{'}\setminus \{x\} ]&=&
 \mu[ \xi(x) \, \, \mid \, \, \xi(.) \, \hbox{on} \, \, \{y_1, y_2, \cdots, y_d\}\cup \bigl(\bigcap_{i=1}^{d} T_{y_i}^{'}(x) \cap V^{'}\bigr)  ]\\
&=&  \mu[ \xi(x) \, \, \mid \, \, \xi(.) \, \hbox{on} \, \, N_x\cap V^{'} ]
\end{eqnarray*}
 because
\begin{equation}\label{emptyset}
 \bigcap_{i=1}^{d} T_{y_i}^{'}(x)\cap V^{'} = \bigcap_{y\in N_x} T_{y}^{'}(x)\cap V^{'} = \emptyset.
\end{equation}
Therefore, the measure $\mu$ satisfies (\ref{MC_eq}). This finishes the proof, the verification of (\ref{emptyset}) being left to the reader.
\begin{cor}
$$\bigcap_{o\in V}o-\mathcal{BMC}(T) \subseteq \mathcal{MC}(T) \subseteq \mathcal{MF}(T).$$
\end{cor}
\section{One-dimensional BMC}\label{sect_1D}
In this section we consider the one-dimensional lattice $V=\mathbb{Z}$ occupied with its natural structure of tree, where the edge set is $E = \{<k,k+1>, \quad k\in\mathbb{Z}\}$. Here $\Omega = \Xi^{\mathbb Z}$.
\begin{prop}\label{prop_1D}
Let  $\mu$ be a probability measure on $(\Omega, \mathcal{B})$. The following assertions are equivalent:
\begin{description}
  \item{(i)} $\mu$ is a  o-BMC for each $o\in V$;
  \item{(ii)}  $\mu$ is an o'-BMC, for some root $o'\in\mathbb Z$;
  \item{(iii)} $\mu$ is a MC.
\end{description}
In particular, a probability measure on $(\Omega, \mathcal{B})$ is markovian for the backward direction if and  only if it is for the forward direction.
\end{prop}

\textit{Proof.} \\
$(i)\Rightarrow (ii)$ straightforward.\\
$(ii)\Rightarrow (i)$  Let $o'\in\mathbb Z$, without loss of generality we can assume that $o< o'$.
 Observe  that if  $x\ge max(o,o')$ or $x\le min(o,o')$ then  $T_0^{'}(x)= T_{o'}^{'}(x)$. Then (\ref{Mp_x_S_o(x)}) is also true if we replace $o$ by $o'$.\\
Let us  now  examine the case  $o< x <o'$ then $ S^{o'}(x) = \{x-1\}$ and $T_{o'}^{'}(x)= (-\infty, x-1] $. Let $m\in \mathbb N$ and $\xi\in \Omega$.\\
Applying (\ref{Mp_x_S_o(x)}) to $y\ge x$, we get
$$
\mu[\xi(y) \, \, \mid \, \, \xi(y-1), \cdots, \xi(x) ] = \mu[\xi(y) \, \, \mid \,\, \xi(y-1)]
$$
because $\{x, \dots, y-1 \}\subseteq (-\infty, y-1] = \mathbb Z\setminus  T_o^{'}(y)$.\\
According to (\ref{condi_prob_formula}), it follows that
\begin{eqnarray*}
\mu\left[ \xi(x-1)\, \, \mid \,\, \xi(.) \,\, \hbox{on}\,\,   [x, x+m]\right]& =& \frac{\mu[ \xi(.) \,\, \hbox{on}\,\,   [x-1, x+m]]}{\mu[ \xi(.) \,\, \hbox{on}\,\,   [x, x+m]]}\\
&=&\frac{ \mu[\xi(x-1)]\prod_{k=x}^{x+m}\mu[ \xi(k) \,\, \mid \,\,  \xi(k-1)]}{\mu[\xi(x)]\prod_{k=x+1}^{x+m}\mu[ \xi(k) \,\, \mid \,\,  \xi(k-1) ]}\\
&=&\frac{\mu[ \xi(x) \,\, \mid \,\,  \xi(x-1)]\times\mu[\xi(x-1)]}{\mu[\xi(x)]}\\
&=& \mu[\xi(x-1) \, \, \mid \, \, \xi(x)].
\end{eqnarray*}
Thus
$$
\mu[\xi(x) \, \, \mid \,  \xi(.) \, \, \hbox{on} \,\,  \mathbb Z\setminus T_{o'}^{'}(x) ] =  \mu[\xi(x-1) \, \, \mid \, \, \xi(x) ]
$$
 for all $x\in \mathbb Z$. Hence $\mu$ is a $o^{'}-BMC$.
 $(ii) \Rightarrow (iii)$
Let $x\in\mathbb Z$ and $m\in \mathbb N$.   Since  $\mu$ is $BMC$ then it is a $o-BMC $ for $o= x-1$ and $T^{'}(x) = [x+1, \infty)$. By (\ref{Mp_x_S_o(x)}), it follows that
$$
\mu [\, \xi(x) \, \, \mid \, \, \xi(x-1), \cdots, \xi(x-m) ] = \mu[\xi(x) \, \, \mid \, \, \xi(x-1)\,  ].
$$
 Therefore, $\mu$ is a Markov chain.\\
$(iii) \Rightarrow (i)$ If $\mu$ is a Markov chain then
$$
\mu[\xi(x) \, \, \mid \, \, \xi(.) \, \, \hbox{on} \, \, (-\infty, x-1)].
$$
By taking $x>0$, this implies that $\mu$ is $0$-block Markov chain, which completes the proof.
\begin{rem}
   Proposition \ref{prop_1D} may be summarized by saying that for each $x\in\mathbb Z$ the triple $(\mathcal F_{[x+1, \infty)}, \mathcal F_{x}, \mathcal F_{(-\infty, x-1]})$ is a Markov triple in the sense of (\ref{Mt}). Namely, this result is still true by taking  $\mathbb{N}$ instead of $\mathbb{Z}$. However, a slight modification on the one dimensional lattice can provide a counter-example in  the multi-dimensional case, in fact we have the following section.
 \end{rem}
\section{Counter-example}\label{sect_Cexple}
 Consider the sets $V = \mathbb N \times \{0\} \cup \{(0,1) , (0, -1)\}\subset \mathbb{Z}^2$ and $E= \{ \{x,y\} \in V \, ; \, |x-y| =1 \}$ where $|(a,b)| = |a| + |b|$.
We get then The tree $T= (V, E)$   (see Fig.\ref{Tgraph}).
 \begin{center}
 \includegraphics[width=0.9\textwidth]{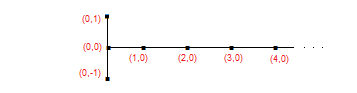}
\end{center}
Consider a  $\{0,1\}$-valued  Markov chain  $(X_n)_{n\ge 0}$ with initial measure $\mu_0 = \frac{1}{2} (\delta_0 + \delta_1)$ and transition matrix $P = \left[
    \begin{array}{cc}
      1/2 & 1/2 \\
      1 & 0 \\
    \end{array}
  \right]$.

Define the $\{0,1\}$-valued  stochastic process $(Z_{u})_{u\in V}$ by
$$
Z_u = \left\{
   \begin{array}{lll}
     X_0, & \hbox{if} &u= (0,-1) ; \\
     X_1, & \hbox{if } &  u= (0,0) ; \\
      X_2, & \hbox{if} &  u= (0,1) ;\\
     X_{n+2}, & \hbox{if} &  u= (n,0),\, n\ge 1.
   \end{array}
 \right.
$$
 Let $\Xi= \{0,1\}$ and $\mu$ be the probability measure on $\Xi^{V}$ associated with $(Z_u)_{u\in V}$. Let $o= (0,-1)$ and $o' = (0,1)$, it easy to check that $\mu$ is an $o$-BMC. However, $\mu$ is not a $o'$-BMC. In fact, if $x= (0,0)$ we have $S_{o'}(w) = \{ (0, -1), (1,0)\},\, r(x)= (0, 1)$   and $T_{o'}^{'}(x) = V\setminus\{x , r(x)\} $ . Let $\xi \equiv 0 \in \Omega$
$$
 \mu[\xi(.) \,\, \hbox{on} \, S_{o'}(x) \, \mid \, \xi(.) \,\, \hbox{on} \, V\setminus T^{'}_{o'}(x)]
$$
$$
 =  \mathbb P[Z_{(0, -1)} =0,\,  Z_{(1,0)} =0 \, \mid \, Z_{(0, 0)} = 0,\,  Z_{(0,1)} =0 ] =  \frac{1}{6}.
$$
On the other hand
$$
\mu[\xi(.) \, \, \hbox{on}\, \, S_{o'}(x) \, \, \mid \, \, \xi(x) ] = \mathbb P[Z_{(0,-1)}=0, Z_{(1,0)} = 0 \, \, \mid  \, \, Z_{(0,0)} =0  ]
= \frac{1}{4}.
$$
This leads to
$$
 \mu[\xi(.) \,\, \hbox{on} \, S_{o'}(x) \, \mid \, \xi(.) \,\, \hbox{on} \, V- T^{'}_{o'}(x) ] \ne \mu[\xi(.) \, \, \hbox{on} S_{o'}(x) \, \, \mid \, \, \xi(x) ].
$$
 Hence $\mu$ is not an $o'$-BMC.\\
Furthermore, the probability measure $\mu$ is not a MC. In fact, by  considering  the subtree with vertex set $V_0  = \{(0,1), (0,0), (1,0) \}$.  We get
$$
\mu[\xi((1,0)) \, \mid \,   \xi((0,0)), \xi((0,1)) ] = \frac{1}{2} \ne \frac{3}{4} = \mu[\xi((1,o)) \, \mid \,   \xi((0,0)) ].
$$

%% file: osid_BMC17082020.bbl
\begin{thebibliography}{91}
\bibliographystyle{plain}
\bibitem{MS19} Mukhamedov F., Souissi A., Quantum Markov States on Cayley trees, \textit{J. Math. Anal. Appl.}{\bf 473} (2019) 313-333.
\bibitem{MS20}\emph{ Mukhamedov F., Souissi A.}, Diagonalizability of quantum Markov States on trees, to appear \textit{J. stat. phys.} (2020).

 \bibitem{HiWeZhi18} H. Huilin, Y.  Weiguo,  S. Zhiyan  \textit{The Shannon–McMillan theorem for Markov chains in Markovian environments indexed by homogeneous trees}, Communications in Statistics - Theory and Methods 47:21, 5286-5297,(2018)
\bibitem{D68} R.L. Dobrushin, \textit{The Description of a Random Field by Means of Conditional Probabilities and Conditions of Its Regularity}. Probab. Theory Appl. 13, 201-22 (1968).
\bibitem{Spa75} A. Spataru, \textit{Construction of a Markov Field on an infinite tree},  Adv. Math. 81, 105-116 (1990).
\bibitem{Spi75} F. Spitzer,  \textit{Markov random fields on an infinite tree }, Ann. Probab.3, 387-398 (1987) .
\bibitem{Za83} S. Zachary, \textit{Countable state space Markov random fields and Markov chains on trees}, Ann. Probab. 4, 894-903 (1983).
\bibitem{Geor} Georgi H.-O. \textit{Gibbs measures and phase transitions}, de Gruyter Studies in Mathematics vol. 9, Walter de Gruyter, Berlin, 1988.

\bibitem{ChGMMP2013} N.Chandgotia, G. Han, B. Marcus, T. Meyerovitch, R. Pavlov \textit{One-dimensional Markov random fields, Markov chains and topological Markov fields},  Am. Math. S  142,  1, 227-242 (2014).
\bibitem{Lau98} S.L. Lauritzen,  \textit{Graphical Models}, Oxford university press (1996)
\bibitem{Ru_Ly16} R. Lyons, Y. Peres,  \textit{Probability on Trees and Networks}, Combridge university press (2016).
\bibitem{Minlos2014} R. A. Minlos, E. A. Pecherskii, S. A. Pirogov, \textit{Gibbs Random Fields on a Lattice: Definitions, Existence,
Uniqueness, and Phase Transitions},  J. Comm. Tech. and Elec., 59, 6, 576-594 (2014).
\bibitem{Rozikov} U. A. Rozikov,  \textit{Gibbs Mesures on Cayley trees}, World scientific (2013)
\bibitem{MuAk13} Mukhamedov, F.,  Akın, H., \textit{ Phase transitions for p-adic Potts model on the Cayley tree of order three}. Journal of Statistical Mechanics: Theory and Experiment, 2013(07), 30
   \bibitem{MuSaKh15} Mukhamedov, F., Saburov, M.,   Khakimov, O., \textit{ On p-adic Ising–Vannimenus model on an arbitrary order Cayley tree}. Journal of Statistical Mechanics: Theory and Experiment, 2015(5), 26
\bibitem{Coa} J. Cao, K. J. Worsley \textit{Applications of Random Fields in Human Brain Mapping}, Spatial Statistics: Methodological Aspects and Applications. Lecture Notes in Statistics, vol 159. Springer, New York, NY (2001).
    \bibitem{Gly} C. Glymour, \textit{The Mind's Arrows: Bayes Nets and Graphical Causal Models in Psychology},   The MIT Press (2001).
 \bibitem{Norris}     J. R. Norris,  \textit{Markov chains, Cambridge Series in Statistical and Probabilistic Mathematics,}
vol. 2, Cambridge University Press, Cambridge, 1998. Reprint of 1997 original. MR1600720 (99c:60144).
\end{thebibliography}
